\documentclass[aps,prl,amsmath,amssymb,onecolumn,nofootinbib,showpacs,tightenlines,14pt]{revtex4}

\usepackage{epsfig,graphicx}
\usepackage{bm}
\newcommand{\be}{\begin{equation}}
\newcommand{\ee}{\end{equation}}
\newcommand{\bea}{\begin{eqnarray}}
\newcommand{\eea}{\end{eqnarray}}

\begin{document}
\title{Semi periodic maps on complex manifolds}
\author{Ali Reza Khatoon Abadi}
 \affiliation{Department of
Mathematics, Islamic azad university,Tehran's West
branch,Tehran,IRAN}
\author{H.R.Rezazadeh}\author{F.Golgoii}
\email{h-rezazadeh@kiau.ac.ir} \affiliation{Department of
Mathematics,Islamic Azad university,KARAJ branch,Karaj,IRAN}

\begin{abstract}
In this letter we proved this theorem: \emph{if $F$ be a holomorphic
mapping of $T_{\Omega}$ to a mapping manifold $X$ such that for
every compact subset $K\subset \Omega$ the mapping $F$ is uniformly
continues on $T_{K}$ and $F(T_{K})$ is a relatively compact subset
of $X$. If the restriction of $F(z)$ to some hyperplane
$\mathbb{R}^{m}+iy'$ is semi periodic, then $F(z)$ is an semi
mapping of $T_{\Omega}$ to $X$.}
\end{abstract} \maketitle

\section{History}
In mathematics, an almost(semi) periodic function is, loosely
speaking, a function of a real number that is periodic to within any
desired level of accuracy, given suitably long, well-distributed
"almost-periods". The concept was first studied by Harald Bohr and
later generalized by Vyacheslav Stepanov, Hermann Weyl and Abram
Samoilovitch Besicovitch, amongst others. There is also a notion of
almost periodic functions on locally compact abelian groups, first
studied by John von Neumann.

Almost periodicity is a property of dynamical systems that appear to
retrace their paths through phase space, but not exactly. An example
would be a planetary system, with planets in orbits moving with
periods that are not commensurable (i.e., with a period vector that
is not proportional to a vector of integers). A theorem of Kronecker
from diophantine approximation can be used to show that any
particular configuration that occurs once, will recur to within any
specified accuracy: if we wait long enough we can observe the
planets all return to within a second of arc to the positions they
once were in.There are several different inequivalent definitions of
almost periodic functions. An almost periodic function is a
complex-valued function of a real variable that has the properties
expected of a function on a phase space describing the time
evolution of such a system. There have in fact been a number of
definitions given, beginning with that of Harald Bohr. His interest
was initially in finite Dirichlet series. In fact by truncating the
series for the Riemann zeta function $\zeta(s)$ to make it finite,
one gets finite sums of terms of the type
\begin{eqnarray}
    e^{(\sigma+it)\log n}\,
    \end{eqnarray}
    with s written as $\sigma+it$ – the sum of its real part ? and imaginary part it. Fixing $\sigma$
     so restricting attention to a single vertical line in the complex plane, we can see this also as
\begin{eqnarray}
    n^\sigma e^{(\log n)it}.\,
\end{eqnarray}
Taking a finite sum of such terms avoids difficulties of analytic
continuation to the region$\sigma < 1$. Here the 'frequencies' log n
will not all be commensurable (they are as linearly independent over
the rational numbers as the integers n are multiplicatively
independent – which comes down to their prime factorizations).

With this initial motivation to consider types of trigonometric
polynomial with independent frequencies, mathematical analysis was
applied to discuss the closure of this set of basic functions, in
various norms.

The theory was developed using other norms by Besicovitch, Stepanov,
Weyl, von Neumann, Turing, Bochner and others in the 1920s and
1930s.

\section{Definitions and some theorems}
A contious mapping $F$ of a tube
\begin{eqnarray} T_{K}={z=x+iy:x\in
\mathbb{R}^{m},y\in K\subset \mathbb{R}^{m}}
\end{eqnarray}
to a metric space $X$ is semi periodic if the family ${F(z+t)}_{t\in
\mathbb{R}^{m}}$ of shifts along $\mathbb{R}^{m}$ is a relatively
compact set with respect to the topology of the uniform
convergence of $T_{K}$.\\
Further, let X be a manifold , and F be a holomorphic mapping of a
tube
\begin{eqnarray}
T_{\Omega}={z=x+iy:x\in \mathbb{R}^{m},y \in \Omega}
\end{eqnarray}
with the convex open base $\Omega\subset \mathbb{R}^{m}$, to X. We
will say that F is semi periodic if the restriction of F to each
tube $T_{K}$ , with the compact base $K\subset\Omega$ is semi
periodic.\\
For $X=C$ we obtain the well-known class of holomorphic semi
periodic functions; for $X=C^{q}$ the corresponding class was being
studied in \cite{Favorov,Raskovskii1,Raskovskii2}; for $X=CP$ we get
the class of meromorphic semi periodic functions that was being
studied in \cite{Sunyer}; the class of holomorphic semi periodic
curves, corresponding to the case $X=CP^{q}$, was being studied in
\cite{Parfyonova}.
\subsection{\textbf{Uniform or Bohr or Bochner almost periodic functions}}
The following theorem is due to Bohr\cite{Bohr}:
\\
: Bohr\cite{Bohr}
 defined the uniformly almost-periodic functions as the
closure of the trigonometric polynomials with respect to the uniform
norm
\begin{eqnarray}
    ||f||_\infty = \sup_x|f(x)|
\end{eqnarray}
(on continuous functions f on R). He proved that this definition was
equivalent to the existence of a relatively-dense set of $\epsilon$
almost-periods, for all $\epsilon > 0$: that is, translations
$T(\epsilon) = T$ of the variable t making
\begin{eqnarray}
    \left|f(t+T)-f(t)\right|<\varepsilon.
\end{eqnarray}

An alternative definition due to Bochner (1926) is equivalent to
that of Bohr and is relatively simple to state:

    A function f is almost periodic if every sequence ${ƒ(t_{n} + T)}$
     of translations of f has a subsequence that converges uniformly for T in $(-\infty,\infty)$.

The Bohr almost periodic functions are essentially the same as
continuous functions on the Bohr compactification of the reals.
 Stepanov almost periodic functions

The space Sp of Stepanov almost periodic functions  was introduced
by V.V. Stepanov \cite{Stepanov}. It contains the space of Bohr
almost periodic functions. It is the closure of the trigonometric
polynomials under the norm
\begin{eqnarray}
    ||f||_{S,r,p}=\sup_x \left({1\over r}\int_x^{x+r} |f(s)|^p \, ds\right)^{1/p}
\end{eqnarray}
for any fixed positive value of r; for different values of r these
norms give the same topology and so the same space of almost
periodic functions (though the norm on this space depends on the
choice of r).  Weyl almost periodic functions

The space Wp of Weyl almost periodic functions  was introduced by
Weyl . It contains the space Sp of Stepanov almost periodic
functions. It is the closure of the trigonometric polynomials under
the seminorm
\begin{eqnarray}
    ||f||_{W,p}=\lim_{r\mapsto\infty}||f||_{S,r,p}
\end{eqnarray}
Warning: there are nonzero functions ƒ with $||ƒ||W,p = 0$, such as
any bounded function of compact support, so to get a Banach space
one has to quotient out by these functions.  Besicovitch almost
periodic functions

The space Bp of Besicovitch almost periodic functions was introduced
by Besicovitch (1926). It is the closure of the trigonometric
polynomials under the seminorm
\begin{eqnarray}
    ||f||_{B,p}=\limsup_{x\mapsto\infty}\left({1\over 2x}\int_{-x}^x |f(s)|^p \, ds\right)^{1/p}
\end{eqnarray}
Warning: there are nonzero functions ƒ with ||ƒ||B,p = 0, such as
any bounded function of compact support, so to get a Banach space
one has to quotient out by these functions.

The Besicovitch almost periodic functions in B2 have an expansion
(not necessarily convergent) as
\begin{eqnarray}
    \sum a_ne^{i\lambda_n}
\end{eqnarray}

Conversely every such series is the expansion of some Besicovitch
periodic function (which is not unique).

The space Bp of Besicovitch almost periodic functions

contains the space Wp of Weyl \cite{Weyl}almost periodic functions.
If one quotients out a subspace of "null" functions, it can be
identified with the space of Lp functions on the Bohr
compactification of the reals. Almost periodic functions on a
locally compact abelian group

With these theoretical developments and the advent of abstract
methods (the Peter–Weyl theorem, Pontryagin duality and Banach
algebras) a general theory became possible. The general idea of
almost-periodicity in relation to a locally compact abelian group G
becomes that of a function F in $L^{*}(G)$, such that its translates
by G form a relatively compact set. Equivalently, the space of
almost periodic functions is the norm closure of the finite linear
combinations of characters of G. If G is compact the almost periodic
functions are the same as the continuous functions.

The Bohr compactification of G is the compact abelian group of all
possibly discontinuous characters of the dual group of G, and is a
compact group containing G as a dense subgroup. The space of uniform
almost periodic functions on G can be identified with the space of
all continuous functions on the Bohr compactification of G. More
generally the Bohr compactification can be defined for any
topological group G, and the spaces of continuous or Lp functions on
the Bohr compactification can be considered as almost periodic
functions on G. For locally compact connected groups G the map from
G to its Bohr compactification is injective if and only if G is a
central extension of a compact group, or equivalently the product of
a compact group and a finite-dimensional vector space.

\emph{If a holomorphic bounded function on a strip is semi periodic
on some straight line in this strip, then this function is semi
periodic on the whole strip}.\\
This theorem was extended to holomorphic functions on a tube
in\cite{Udodova}; besides usual uniform metric, various integral
metric were being studied here.\\
The direct generalization of Bohr's theorem to complex manifold is
not valid.
\section{Theorem}
\textbf{Theorem:}\emph{if $F$ be a holomorphic mapping of
$T_{\Omega}$ to a mapping manifold $X$ such that for every compact
subset $K\subset \Omega$ the mapping $F$ is uniformly continues on
$T_{K}$ and $F(T_{K})$ is a relatively compact subset of $X$. If the
restriction of $F(z)$ to some hyperplane $\mathbb{R}^{m}+iy'$ is
semi periodic, then $F(z)$
is an semi mapping of $T_{\Omega}$ to $X$.}\\
Also there is a corollary:\\
\textbf{Corollary}:\emph{Let F be a holomorphic mapping form
$T_{\Omega}$ to a compact complex manifold X such that F is
uniformly continues on $T_{K}$ for every compact set $K\subset
\Omega$. If the restriction of F(z) to some hyperplane
$\mathbb{R}^{m}+iy'$ is semi periodic , then F(z)
is an semi periodic mapping of $T_{\Omega}$ to X.}\\
 \textbf{Proof of the theorem:}\\
 Take an arbitrary sequence $\{t_{n}\}\subset \mathbb{R}^{m}$ .
 Since the function F(z) is uniformly continues, the family
 $\{F(z+t_{n})\}$ is equicontinous on each compact set $S\subset
 T_{\Omega}$. Further, it follows from the condition of the Theorem
 that the union of all the images of S under mappings of this family
 is contained compact subset of X. Therefore, passing on to a
 subsequence if necessary, we may assume that the sequence
 $\{F(z+t_{n})\}$ converges  to a holomorphic mapping G(z) uniformly
 on every compact subset of $T_{\Omega}$. It easy to see that the
 mapping G(z) is bounded and uniformly continues on every tube
 $T_{K}$ with the compact base $K\subset\Omega$. Let us prove that
 this convergence is uniform on every $T_{K}$. Assume the contrary.
 Then we get
 \begin{eqnarray}
 ||F(z+t_{n}),G(z_{n})||>0
\end{eqnarray}
for some sequence $z_{n}=x_{n}+iy_{n}\in T_{K'}$, where $K'$ is some
compact subset of $\Omega$.\\
Replacing sequences by a subsequence if necessary, we may assume
that the mapping $G(x_{n}+z)$ converge to a holomorphic mapping
$H(z)$, and the mappings $F(z+x_{n}+t_{n})$ converge to a
holomorphic mapping $\tilde{H}(z)$ uniformly on every compact
subsets of $T_{\Omega}$. We may also assume that $y_{n}\rightarrow
y_{0}\in K'$. Using (1) we get
\begin{eqnarray}
|\tilde{H}(iy_{0})-H(iy_{0}|\geq\varepsilon_{0}
\end{eqnarray}
Since the mapping $F(x+iy')$ of $\mathbb{R}^{m}$ to X is semi
periodic, we may assume that a subsequence of mappings
$F(x+t_{n}+iy')$ converges to $G(x+iy')$ uniformly in
$x\in\mathbb{R}^{m}$. Therefor the sequences of mappings
$F(x+x_{n}+t_{n}+iy')$ and $G(x+x_{n}+iy')$ have the same limit,
i.e.,$\tilde{H}(x+iy')=H(x+iy')$ for all $x\in \mathbb{R}^{m}$.
Since $\tilde{H}(z),H(z)$ are holomorphic mappings we get
$\tilde{H}(x+iy')\equiv H(x+iy')$ on $T_{Omega}$. This contradiction
proves the Theorem.

\end{document}